\documentstyle[amsfonts,11pt]{article}

\def\Iso {\mathop{\hbox{Iso}}}

\newtheorem{teor}{Theorem}[section]
\newtheorem{defin}[teor]{Definition}
\newtheorem{remar}[teor]{Remark}

\newtheorem{corol}[teor]{Corollary}
\newtheorem{lemma}[teor]{Lemma}

\newtheorem{Example}[teor]{Example}

\newcommand{\fdim}{\hspace*{\fill}$\Box$}
\newcommand{\dimostr}{{\bf Proof: }}

\newcommand{\complex}{\Bbb{C}}

\newcommand{\K}{K\"{a}hler}

\begin{document}

\noindent
\centerline {\LARGE\bf \K\ manifolds
and their relatives}





\vspace{0.5cm}

\centerline{\small Antonio J. Di Scala}
\centerline{\small Dipartimento di Matematica
-- Politecnico di Torino -- Italy}
\centerline{\small e-mail address: antonio.discala@polito.it}

\vspace{0.3cm}
\centerline{\small and}
\vspace{0.3cm}

\centerline{\small Andrea Loi}
\centerline{\small Dipartimento di Matematica e Informatica
-- Universit\`{a} di Cagliari -- Italy}
\centerline{\small e-mail address: loi@unica.it}

\vspace{0.3cm}

\begin{abstract}
\noindent
Let $M_1$ and $M_2$ be two \K\ manifolds.
We call $M_1$ and $M_2$
{\em relatives} if they share a non-trivial \K\ submanifold
$S$, namely, if there exist two holomorphic and isometric immersions
(\K\ immersions) $h_1: S\rightarrow M_1$ and $h_2: S\rightarrow
M_2$. Moreover, two \K\ manifolds $M_1$ and $M_2$ are said to be {\em
weakly relatives} if there exist two locally isometric (not
necessarily holomorphic) \K\ manifolds  $S_1$ and $S_2$   which
admit two \K\ immersions into $M_1$ and $M_2$ respectively.
The notions introduced are not equivalent (cfr. Example \ref{hyperk}).
Our main results in this paper  are Theorem \ref{mainteor} and
Theorem \ref{mainprop}.
In the first theorem we show that a complex bounded domain $D\subset {\complex}^n$
with its Bergman metric
and a  projective \K\  manifold
(i.e. a projective manifold  endowed with the restriction of the  Fubini--Study metric)
are not relatives.
In the second theorem we prove that a Hermitian symmetric space of
noncompact type  and a projective \K\ manifold  are not weakly relatives.
Notice that the proof of the second result
does not follows trivially from the first one.
We also remark that the above results are of local nature,
i.e. no assumptions are used about the compactness or completeness of the manifolds involved.

\vspace{0.3cm}


\noindent
{\it{Subj.Class}}: 53C55, 58C25.

\end{abstract}
\section{Introduction}

The study of holomorphic and isometric immersions
between \K\ manifolds (called \K\ immersion in the sequel)
was started by Eugenio Calabi who,
in his pioneering work \cite{ca} of $1953$,
solved the problem of deciding
about the existence of \K\ immersions
between complex space forms.
More specifically, he proved that two
complex space forms with curvature  of different
sign cannot be \K\ immersed one into another
and, in particular that,
for complex space forms of  the same type,
just projective spaces can be
embedded between themselves in a non trivial way
by using Veronese mappings.
Unfortunately, this subject has not been
further explored by other authors
as pointed out by Marcel Berger in \cite{berger}
who referring to Calabi's paper wrote:
{\em \lq\lq this wonderful text remains quite unknown
and almost unused..."}.

The authors believe that the results of the present
paper are in accordance with Berger's opinion
(see also page $528$ in \cite{bepan})
and could stimulate
future research in this field.

\vskip 0.3cm

\noindent
In order to state our first result
we give the following:

\begin{defin}
Let $r$ be a positive integer.
Two \K\ manifolds $M_1$ and $M_2$ are said to be {\em $r$-
relatives} if they have in common a complex $r$-dimensional \K\ submanifold
$S$, i.e. there exist two \K\  immersions $h_1:S\rightarrow M_1$ and
$h_2:S\rightarrow M_2$.
Otherwise, we say that $M_1$ and $M_2$ are
not relatives.
\end{defin}

Our first result is the following:

\begin{teor}\label{mainteor}
A bounded domain $D\subset {\complex}^n$
endowed with its Bergman metric
and a  projective \K\  manifold endowed with the restriction of  the Fubini--Study metric are not relatives.
\end{teor}

The following defintion  generalizes the previous one:

\begin{defin}
Let $r$ be a positive integer.
Two \K\ manifolds $M_1$ and $M_2$
are said to be {\em weakly $r$-relatives} if there exist
two locally isometric \K\ manifolds  $S_1$ and $S_2$ of complex dimension $r$  which admit two \K\ immersions
$h_1:S_1\rightarrow M_1$ and $h_2:S_2\rightarrow M_2$.
\end{defin}

\noindent
We remark that here the local isometry
between $S_1$ and $S_2$ is not assumed to be holomorphic (cfr.
Example \ref{hyperk} below).

\vskip 0.3cm

\noindent
Let us now state our second  result.

\begin{teor}\label{mainprop}
An irreducible Hermitian symmetric space
of noncompact type
and a projective \K\  manifold   are not
weakly relatives.
\end{teor}

Notice that even if the word
{\em weakly
relatives} in the above theorem
is replaced by the word  {\em relatives}
the statement  does not directly follow from Theorem \ref{mainteor}.
The reason is that
the metric on  Hermitian symmetric spaces
of noncompact type  is either
the Bergman metric or a non-trivial multiple of it.

\vskip 0.3cm

Since an irreducible Hermitian symmetric space of compact
type admits a \K\
embedding into a complex projective space
(see e.g. \cite{diastherm} and
\cite{kobounded}),
Theorem
\ref{mainprop} yield  to
the following appealing:

\begin{corol}
An irreducible Hermitian symmetric space of
noncompact type and an irreducible
Hermitian symmetric space of compact type,
are not weakly relatives.
\end{corol}

\section{Proof of the main results}
Let $\Phi_B (z,  z)$ be
a (globally defined)
\K\ potential for the Bergman metric  $g_B$ on a bounded
domain $D\subset {\complex}^n$.
Then $\Phi_B (z,  z)=\log K_B(z, z)$
where $K_B(z,  z)$ is the Bergman kernel function
on $D$, namely
$K_B(z,  z)=\sum_{j=0}^{+\infty}|F_j(z)|^2$
where
$F_j, j=0,1,\dots $ is  an orthonormal basis
for the Hilbert space ${\cal H}$
consisting of  square integrable
holomorphic functions on $D$.
Observe that, from the boundedness
of $D$, ${\cal H}$
contains all polynomials.
In particular each element
of the sequence $z_1^k, k=0,1,\dots$
belongs to ${\cal H}$, where $z_1$
is the first variable of $z=(z_1,\dots , z_n)\in D\subset {\complex}^n$.
By applying the Gram--Schmidt
orthonormalization procedure to the sequence
$z_1^k$ we can assume that
there exists a
sequence of linearly independent  polynomials
$P_k(z_1), k=0,1\dots$
in the variable $z_1$
such that $P_0(z_1)=F_0(z_1,\dots ,z_n)=\lambda_0\in {\complex}^*$
and $P_k(z_1)=F_k(z_1,\dots ,z_n), \forall k=1,\dots$.
Consider now the holomorphic map
of $D$ into the standard complex Hilbert space
$l^2({\complex})$
given by:

\begin{equation}\label{pphi}
\phi :D\rightarrow l^2({\complex}),
z=(z_1,\dots , z_n)\mapsto
(P(z_1), F(z)),
\end{equation}
where
$P(z_1)=(P_0(z_1),P_1(z_1),\dots )$
and
$F(z)$ is the sequence obtained by
deleting
the sequence
$z_1^k$
from the sequence $F_j(z)$.

Observe that $l^2({\complex})$
can be seen   as the affine chart $Z_0\neq 0$ of the
infinite dimensional
complex projective space
${\complex}P^{\infty}$,
endowed with homogeneous coordinates
$[Z_0,\dots , Z_j,\dots ]$.
Moreover, the Fubini-Study metric $g_{FS}^{\infty}$
of ${\complex}P^{\infty}$
restricts to the \K\ metric
$$\frac{i}{2}\partial\bar
\partial\log (1+\sum_{j=1}^{+\infty}|w_j|^2)\
(w_j=\frac{Z_j}{Z_0})$$
on $l^2({\complex})$ and
it follows by the very definition of the Bergman metric
that the mapÊ
(\ref{pphi})
induces a  \K\ immersion
\begin{equation}\label{phi}
\Phi(z)=[P(z_1), F(z)]:(D, g_B)\rightarrow
({\complex}P^{\infty}, g_{FS}^{\infty}).
\end{equation}

\begin{remar}\rm
The fact that a bounded domain
endowed with its Bergman
metric admits a \K\ immersion
$\Phi$ into the
infinite dimensional
complex projective space is
well-known and was first pointed
out by Kobayashi \cite{kobounded}.
For the proof of our main result  it is crucial
that the map $\Phi$ can be put in the special form
(\ref{phi}).
\end{remar}

For later use we give the following definition.
Let $S$ be a complex manifold.
We say that a holomorphic map
$\Psi :S\rightarrow {\complex}P^{\infty}$
is {\em non-degenerate}
iff $\Psi (S)$ is not contained in any
{\em finite} dimensional
complex projective space
${\complex}P^N\subset {\complex}P^{\infty}$.
The following lemma summarizes what we need
about non-degenerate maps.

\begin{lemma}\label{lemmafull}
Let $S\subset {\complex}^n$
be an open subset of ${\complex}^n$
and let
$$\Psi:S\rightarrow {\complex}P^{\infty}:
z\mapsto [\psi_0(z) ,\psi_1(z) ,\dots]$$
be a
holomorphic map
induced by  the holomorphic map
$$\psi:S\rightarrow l^2({\complex}):
z\mapsto (\psi_0(z) ,\psi_1(z) ,\dots)$$
where $\psi_j, j=0, 1\dots$
is an infinite sequence
of holomorphic functions on $S$.
Assume that there exists an infinite
subsequence $\psi_{j_{\alpha}}$ of $\psi_j$,
consisting of  linearly
independent functions
such that for all $s\in S$
there exists a function of this subsequence non-vanishing at $s$.
Then
$\Psi$ is non-degenerate.
Furthermore,
if $\Psi$ is non-degenerate and
$\tilde\Psi:S\rightarrow {\complex}P^{\infty}$
is another holomorphic immersion
which induces on $S$
the same \K\ metric induced by $\Psi$,
i.e. $\Psi^*(g_{FS})=\tilde\Psi^*(g_{FS})$,
then also
$\tilde\Psi$ is non-degenerate.
\end{lemma}
\dimostr
Let $W$ be the infinite dimensional
complex subspace of $l^2({\complex})$ spanned by the vectors
$e_{j_{\alpha}}$, where $e_j$ is the canonical basis of $l^2({\complex})$. Denote by $\pi :l^2({\complex})\rightarrow W$
the projection onto $W$. Therefore, the map
$\pi\circ \psi:S\rightarrow W\subset l^2({\complex})$
induces a holomorphic and non-degenerate map $S\rightarrow {\complex}P^{\infty}$. Hence, {\em a fortiori},
the map $\Psi$ is non-degenerate.
The proof of the  second part of the lemma is an
immediate consequence of Calabi's rigidity
theorem (see \cite{ca}), which asserts that any two \K\ immersions
$\Psi_1, \Psi_2$ of a \K\ manifold $S$ into ${\complex}P^{\infty}$
are related by a unitary transformation $U$ of ${\complex}P^{\infty}$,
i.e. $U\circ\Psi_1=\Psi_2$.
\fdim

\vskip 0.3cm

\noindent
{\bf Proof of Theorem \ref{mainteor}}

\noindent
We can restrict ourself to prove that
the domain
$D\subset {\complex}^n$ equipped with its Bergman
metric is not relative
to any complex projective space
${\complex}P^m$.
Assume by contradiction this is the case.
Then,  there
exists an open subset $S\subset {\complex}$
passing through the origin
and two  \K\ immersions
$f:S\rightarrow D$ and
$h:S\rightarrow {\complex}P^m$.
If $(f_1, \dots, f_n)$
denote the components of $f$
we can  assume that
$\frac{\partial f_1}{\partial \xi}(0)\neq 0$,
where $\xi$ is the complex coordinate on $S$.
Consider the
\K\  immersion of $S$ into
${\complex}P^{\infty}$ given by the composition
$\Phi\circ  f:S\rightarrow {\complex}P^{\infty}$,
where $\Phi$ is given by (\ref{phi}).
We claim that this map is indeed non-degenerate.
In order to prove our claim
observe that from (\ref{phi}) one gets:
$$(\Phi\circ  f)(\xi)=[P(f_1(\xi)), F(f_1(\xi),\dots ,f_n(\xi))].$$
Hence, by the first part of Lemma \ref{lemmafull}, it is enough to prove that $P_k(f_1(\cdot)), k=0, 1,\dots$
is a sequence of linearly independent functions
on $S$.
So, let $q$ be any positive integer and assume
that there exist $q$ complex numbers $a_0,\dots ,a_q$
such that
\begin{equation}\label{eqa}
a_0P_0(f_1(\xi))+\cdots a_qP_q(f_1(\xi))=0,\ \forall\xi\in S.
\end{equation}
By the assumption on $f_1:S\rightarrow {\complex}$ it follows
that $f_1 (S)$ is on open subset of ${\complex}$.
Therefore, equality (\ref{eqa}) is satisfied on all ${\complex}$,
and since $P_1,\dots P_q$ are linearly independent all the $a_j$'s
are forced to be zero, proving our claim.
Next, consider the \K\ immersion  of
$S$ into
${\complex}P^{\infty}$ given by the composition
$i\circ h$,
where  $i:{\complex}P^m \hookrightarrow {\complex}P^{\infty}$
is the natural inclusion.
Since this map
is obviously degenerate
the second part of Lemma \ref{lemmafull} yields the desired contradiction.
\fdim

\vskip 0.3cm

\noindent
Before proving Theorem \ref{mainprop} let us explain the two main problems
one has to face in proving it.

\vskip 0.3cm
The first one comes from the fact that
a Hermitian symmetric space of noncompact type
$(D, g)$ is equivalent to a
bounded symmetric domain with its Bergman metric
$(D, g_B)$
up to homotheties, i.e.
$g=cg_B, c>0$.
Hence, we cannot apply  directly Theorem
\ref{mainteor} even when  {\em weakly
relatives} is replaced by  {\em relatives}.
Indeed, one can easily exibit two \K\ manifolds
which are not $r$-relatives (for all $r$) but become $s$-relatives (for some $s$) when one  multiplies
the metric of  one of them by a suitable constant
(for example $({\complex}P^1, g=\lambda g_{FS})$  endowed
with an irrational multiple $\lambda$
of the Fubini--Study metric
is not  $r$-relative to  $({\complex}P^1, g_{FS})$
but the later is obviously  $1$-relative
to itself).

\vskip 0.3cm
The second problem is that weakly relatives \K\
manifolds should not be relatives as shown by the following:

\begin{Example}\rm\label{hyperk}
Let $X$ be a $K3$ surface  with its hyperk\"ahlerian structure
(see e.g. p. $400$ in \cite{besse}). It is well-known that its
isometry group $\Iso (X)$ is finite (see \cite{algr} for a
beautiful description of this group). Let $J_1$ and $J_2$ be two parallel
complex structures which do not belong to the same
$\Iso(X)$-orbit. Then, the two \K\ manifolds $(X, J_1)$ and $(X,
J_2)$ are obviously weakly $2$-relatives but not $2$-relatives.
Moreover, a suitable choice of $J_1$, $J_2$ gives two weakly
relatives \K\ manifolds $(X, J_1)$ and $(X, J_2)$ which are not
relatives i.e. neither $1$-relatives nor $2$-relatives.
\end{Example}

The following lemma allow us to avoid the previous
difficulty.

\begin{lemma}\label{lammarel}
Let $(D, g)$ be a Hermitian symmetric space of noncompact type and
$V$ be a projective \K\ manifold. If $D$ and $V$ are weakly
relatives then $D$ and $V$ are also relatives.
\end{lemma}
\dimostr
Let $L: S_1\rightarrow S_2$ be the local isometry  between
the two \K\ manifolds $S_1$ and $S_2$
which makes $D$ and $V$ weakly relatives and let
$h_1: S_1\rightarrow D$
and $h_2: S_2\rightarrow V$ be the corresponding \K\ immersions.
Since all the concepts involved are of local
nature, we can assume that $L$
is a global isometry and, with the aid of
De Rham decomposition theorem, that the Riemannian manifold $S_1=S_2$ decompose as
$$S_1=S_2=F\times I_1\times\dots  \times I_k,$$
where $F$ is an open subset of the Euclidean space with the flat
metric and $I_j, j=1,\dots k$ are irreducible Riemannian
manifolds. Observe that the factor $F$ is indeed a flat \K\
manifold of $S_2$. Since $S_2$ is a projective \K\ manifold it
follows by the above mentioned theorem of Calabi (see the
introduction) that $F$ is trivial. We also claim that  the above
decomposition does not contain Ricci flat factors. Indeed, assume
for example that $I_j$ is such a factor. Then, as a consequence of
the Gauss equation and the non-positivity of the curvature of $(D,
g)$, it follows that the map $h_1:I_j\rightarrow D$ is totally
geodesic. Since a totally geodesic submanifold of a locally
homogeneous Riemannian manifold is also locally homogeneous, a
well-known theorem of Alekseevsky--Kimel'fel'd--Spiro (see
\cite{alki} and \cite{spiro}) implies that $I_j$ is actually flat,
which proves our claim. Finally, observe that the isometry $L$
above takes an irreducible  factor $I$ of $S_1$ into an irreducible
factor $L(I)$ of $S_2$. Since these factors are not Ricci flat it
is well-known that,  $L : I \rightarrow L(I)$ or its conjugate
$\bar L$ is  holomorphic and so $D$ and $V$ are relatives since
they share the same \K\ manifold $I$. \fdim

\begin{remar}\rm
The above lemma is valid (the proof follows the same line) when
$D$ is a homogeneous bounded domain of non-positive holomorphic
bisectional curvature. Notice that the celebrated example of Pyatetski-Shapiro
\cite{pysh} shows that Hermitian symmetric spaces of non-compact
type are strictly contained in such domains (see also
\cite{datri}).
\end{remar}

\vskip 0.3cm

\noindent
{\bf Proof of Theorem \ref{mainprop}}

\noindent
Let ${\cal H}_m$ be the Hilbert space
consisting of holomorphic functions $f$ on $D$
such that $\int_D\frac{|f|^2}{K_B^{mc}}dz<+\infty$,
where $dz$
is the Lebesgue measure on $D$ and $m$ is a
positive integer.
Let $F_j^m$ be an orthonormal basis
for ${\cal H}_m$. It is not hard to see
that
$\frac{\sum_{j=0}^{+\infty}|F^m_j(z)|^2}
{K^{mc}_B(z, z)}$
is invariant under the action of
the group of isometric
biholomorphisms of $(D, g)$. Since
this group acts transitively on $D$
we have that
\begin{equation}\label{fmB}
\sum_{j=0}^{+\infty}|F^m_j(z)|^2=
c_mK^{mc}_B(z, z), c_m>0.
\end{equation}
It can be shown that
for $m$ sufficiently large,
${\cal H}_m$
contains all polynomials
(see \cite{diastherm} and reference therein).
Fix such an $m$.
As in the proof of the previous theorem
one can built a holomorphic map
\begin{equation}\label{phim}
\Phi_m :D\rightarrow {\complex}P^{\infty},
z=(z_1,\dots , z_n)\mapsto
[P^m(z_1), F^m(z)],
\end{equation}
where $P^m(z_1)=(P^m_0(z_1), P^m_1(z_1),\dots )$
is an infinite sequence of linearly independent
polynomials in the variable $z_1$
and $P^m_0(z_1)$ is a non-zero complex number.
Moreover, it follows by (\ref{fmB}) that
$\Phi_m^*(g_{FS}^{\infty})=mg$.
Observe  now  that
there exists
a holomorphic immersion
$V_m:{\complex}P^m\rightarrow {\complex}P^{\left({N+m\atop N}\right)}$
(obtained by a suitable rescaling of the Veronese embedding)
satisfying $V_m^*(g_{FS})=mg_{FS}$
(see Calabi \cite{ca}).

Assume, by a contradiction, that
the Hermitian symmetric space of noncompact type
$(D, g)$ and ${\complex}P^m$ are weakly relatives.
Then by  Lemma \ref{lammarel}
they are also relatives and so
there exists an open subset $S\subset {\complex}$
and two  \K\ immersions
$f:S\rightarrow D$ and
$h:S\rightarrow {\complex}P^m$.
Then, obviously, $(D, mg)$ and $({\complex}P^m, mg_{FS})$
would be  relatives.
Hence, as in the proof
of Theorem \ref{mainteor},
we get  the desired contradiction,
by applying Lemma \ref{lemmafull}
to the \K\ immersions  $\Phi_m\circ f :S\rightarrow {\complex}P^{\infty}$ and
$i\circ V_m\circ h: S\rightarrow {\complex}P^{\infty}$,
(where $i:{\complex}P^{\left({N+m\atop N}\right)} \hookrightarrow {\complex}P^{\infty}$
is the natural inclusion) which are
respectively non-degenerate and degenerate.
\fdim

\begin{remar}\rm
Observe that Theorem \ref{mainprop} when $D$ is a rank one
Hermitian symmetric space of non-compact type (i.e. $D =
{\complex} H^n$) was prove in \cite{um} by Masaaki Umehara.
Umehara's proof is  based on the
Calabi embedding ${\complex} H^n \hookrightarrow
l^2({\complex})$. Since such (K\"ahlerian) embeddings does not
exist  for higher rank Hermitian symmetric space
of non-compact type (see \cite{dilo} for a proof)
Umehara's approach cannot be used to give an alternative proof
of our theorem.
\end{remar}

\small{}

\end{document}